\documentstyle[twoside]{article}
\title{\bf High-order derivatives of the Bessel functions with an application }
\author{\sc R. B. Paris\footnote{E-mail address:\ \ {\tt r.paris@abertay.ac.uk}}\\
\\
{\em Division of Computing and Mathematics,}\\
{\em Abertay University, Dundee DD1 1HG, UK}\\
}

\begin{document}
\newcommand{\bee}{\begin{equation}}
\newcommand{\ee}{\end{equation}}
\def\f#1#2{\mbox{${\textstyle \frac{#1}{#2}}$}}
\def\dfrac#1#2{\displaystyle{\frac{#1}{#2}}}
\newcommand{\fr}{\frac{1}{2}}
\newcommand{\fs}{\f{1}{2}}
\newcommand{\g}{\Gamma}
\newcommand{\om}{\omega}
\newcommand{\br}{\biggr}
\newcommand{\bl}{\biggl}
\newcommand{\ra}{\rightarrow}
\renewcommand{\topfraction}{0.9}
\renewcommand{\bottomfraction}{0.9}
\renewcommand{\textfraction}{0.05}
\newcommand{\mcol}{\multicolumn}
\newcommand{\gtwid}{\raisebox{-.8ex}{\mbox{$\stackrel{\textstyle >}{\sim}$}}}
\newcommand{\ltwid}{\raisebox{-.8ex}{\mbox{$\stackrel{\textstyle <}{\sim}$}}}
\date{}
\maketitle
\pagestyle{myheadings}
\markboth{\hfill {\it R.B. Paris} \hfill}
{\hfill {\it High-order derivatives of the Bessel functions } \hfill}
\begin{abstract} 
We determine the asymptotic behaviour of the $n$th derivatives of the Bessel functions $J_\nu(a)$ and $K_\nu(a)$, where $a$ is a fixed positive quantity, as $n\to\infty$. These results are applied to the asymptotic evaluation of two incomplete Laplace transforms of these Bessel functions on the interval $[0,a]$ as the transform variable $x\to+\infty$. Similar evaluation of the integrals involving the Bessel functions $Y_\nu(t)$ and $I_\nu(t)$ is briefly mentioned.
\vspace{0.4cm}

\vspace{0.3cm}

\end{abstract}

\vspace{0.2cm}

\noindent $\,$\hrulefill $\,$

\vspace{0.2cm}

\begin{center}
{\bf 1. \  Introduction}
\end{center}
\setcounter{section}{1}
\setcounter{equation}{0}
\renewcommand{\theequation}{\arabic{section}.\arabic{equation}}
In this note we examine the asymptotic behaviour of high derivatives of the Bessel functions $J_\nu(a)$ and $K_\nu(a)$, where $a>0$; that is, we determine order estimates of $J_\nu^{(n)}(a)$ and $K_\nu^{(n)}(a)$ as $n\to\infty$. When $\nu$ is a non-negative integer, we determine an asymptotic expansion for $J_\nu^{(n)}(a)$ as $n\to\infty$. These results are then applied to determine the asymptotic expansion of the integrals
\bee\label{e1}
H_J(a,x)=\int_0^a e^{-xt} J_\nu(t)\,dt \quad(\nu\geq0),\qquad H_K(a,x)=\int_0^a e^{-xt} K_\nu(t)\,dt\quad(0\leq\nu<1)
\ee
as $x\to+\infty$. In the first integral we restrict $\nu$ to satisfy $\nu\geq0$ although it converges for $\nu>-1$; in the second integral we can omit consideration of the case $-1<\nu<0$, since $K_{-\nu}(x)=K_\nu(x)$.

In the application of Watson's lemma \cite[p.~44]{DLMF}, \cite[p.~15]{PHad} to the above integrals for large $x$, the Bessel functions are replaced by their series expansions about $t=0$ and integration is extended to the infinite interval $[0,\infty)$. 
This procedure then simply produces the convergent series expansions with the evaluations
\bee\label{e12}
H_J(a,x)\sim\frac{(\sqrt{1+x^2}-x)^\nu}{\sqrt{1+x^2}}
\ee
and\footnote{When $x=1$, we have the limiting value $H_K(\infty,1)=\pi\nu/\sin \pi\nu$.}
\bee\label{e13}
H_K(a,x)\sim\frac{\pi}{2\sin \pi\nu}\,\frac{1}{\sqrt{x^2-1}}\bl\{(x+\sqrt{x^2-1})^\nu-(x+\sqrt{x^2-1})^{-\nu}\br\}.
\ee
The above values of $H_J(a,x)$ and $H_K(a,x)$ coincide, of course, with the evaluations when the upper limit $a$ is replaced by $+\infty$; see \cite[p.~386(8), 388(9)]{W}.

As an application we employ the asymptotic behaviour of the high-order derivatives of $J_\nu(a)$ and $K_\nu(a)$ to determine the character of the exponentially small contributions to the asymptotic expansion of the integrals in (\ref{e12}) and (\ref{e13}) as $x\to+\infty$. Similar integrals involving the Bessel functions $Y_\nu(t)$ and $I_\nu(t)$ are briefly mentioned.
\vspace{0.6cm}

\begin{center}
{\bf 2. \ The derivatives $J_\nu^{(n)}(a)$ }
\end{center}
\setcounter{section}{2}
\setcounter{equation}{0}
\renewcommand{\theequation}{\arabic{section}.\arabic{equation}}
Let $n$ denote a positive integer and $a>0$ be a parameter. Then the $n$th derivative of the Bessel function $J_\nu(x)$ at $x=a$ can be written by Cauchy's integral formula as
\bee\label{e20}
J_\nu^{(n)}(a)=\frac{n!}{2\pi i} \oint \frac{J_\nu(a+u)}{u^{n+1}}\,du,
\ee
where the integration path is a closed loop surrounding the origin in the positive sense not enclosing the branch point $u=-a$ (when $\nu$ is non-integer). Now expand the contour into a large circular path $C$ of radius $R$, together with paths along the upper and lower sides of the branch cut along $[-a,-\infty)$ connected by a small circular path of radius $\rho$ about the branch point. When $\nu=0, 1, 2, \ldots$ there is no branch point and the integration path can be taken to be the circular circuit of radius $R$.

For $R\to\infty$ we can employ the asymptotic expansion of the Bessel function in the form \cite[p.~228]{DLMF}
\[J_\nu(x)\sim \frac{1}{\sqrt{2\pi x}}\bl\{e^{i\omega}\sum_{k\geq0} \frac{i^k c_k(\nu)}{x^k}+
e^{-i\omega}\sum_{k\geq0}\frac{(-i)^k c_k(\nu)}{x^k}\br\}\]
valid as $|x|\to\infty$ in $|\arg\,x|<\pi$, where $\omega=x-\fs\pi\nu-\f{1}{4}\pi$ and the coefficients $c_k(\nu)$ are
given by
\[c_k(\nu)=\frac{(-)^k(\fs+\nu)_k(\fs-\nu)_k}{2^k k!}.\]
It then follows that for large $|u|$ in $|\arg\,u|<\pi$
\bee\label{e21}
J_\nu(a+u)\sim\frac{1}{\sqrt{2\pi u}}\bl\{e^{i\omega}\sum_{k\geq0} \frac{A_k}{u^{k}}+e^{-i\omega}\sum_{k\geq0} \frac{{\bar A}_k}{u^{k}}\br\},\quad\omega=u+\phi,
\ee
where
\[\phi:=a-\fs\pi\nu-\f{1}{4}\pi,\qquad A_1=ic_1-\fs a_1, \quad A_2=\f{3}{8}a^2-c_2-\f{3}{2}iac_1, \ldots\]
and the bar denotes the complex conjugate.

The contribution ${\cal J}_1$ to (\ref{e20}) from the expanded contour $C$ is then
\[{\cal J}_1\sim\frac{n!}{2\pi i}\,\frac{1}{\sqrt{2\pi}}\bl\{\sum_{k\geq0}A_k\int_C\frac{e^{i\omega}}{u^{k+n+3/2}}\,du+\sum_{k\geq 0}{\bar A}_k\int_C\frac{e^{-i\omega}}{u^{k+n+3/2}}\,du\br\}\]
\bee\label{e22}
=\frac{n! n^{-n-1/2}}{\sqrt{2\pi}}\,\frac{1}{2\pi i}\bl\{e^{i\phi}\sum_{k\geq0}\frac{A_k \Upsilon_k^+(n)}{n^k}+e^{-i\phi}\sum_{k\geq0}\frac{{\bar A}_k \Upsilon_k^-(n)}{n^k}\br\}
\ee
upon making the change of variable $u=nw$.  Here we have introduced
\[\Upsilon_k^\pm(n)=\int_{C'}\frac{e^{n(\pm iw-\log\,w)}}{w^{k+3/2}}\,dw\]
and $C'$ denotes a circular path of radius $R/n$.

With $\psi(w):=iw-\log\,w$, the exponential factor appearing in $\Upsilon_k^+(n)$ has a saddle point (where $\psi'(w)=0$) at $w_s=-i$, where $\psi''(w_s)\equiv \psi_s''=-1$. The path of steepest descent through the saddle is locally horizontal and passes to infinity in the upper half-plane. The path $C'$ can be deformed to pass over the saddle and application of the method of steepest descents then yields
\bee\label{e23}
\Upsilon_k^+(n)\sim e^{n\psi_s} \sqrt{\frac{2\pi}{n}}\,i^{k+3/2} \sum_{r\geq0} \frac{B_r(k)}{n^r}\qquad (n\to\infty),
\ee
where $\psi_s\equiv\psi(w_s)=1+\fs\pi i$ and \cite[pp.~13--14]{PHad}
\[B_0(k)=1,\qquad B_1(k)=-\frac{1}{4\psi_s''}\bl\{2F_2-2\Psi_3F_1+\frac{5}{6}\Psi_3^2-\frac{1}{2}\Psi_4\br\},\]
\[B_2(k)=\frac{3}{16(\psi_s'')^2}\bl\{\frac{2}{3}F_4-\frac{20}{9}\Psi_3F_3+\frac{5}{3}\bl(\frac{7}{3}\Psi_3^2-\Psi_4\br)F_2-\frac{35}{9}\bl(\Psi_3^3-\Psi_3\Psi_4+\frac{6}{35}\Psi_5\br)F_1\]
\[+\frac{35}{9}\bl(\frac{11}{24}\Psi_3^4-\frac{3}{4}\bl(\Psi_3^2-\frac{1}{6}\Psi_4\br)\Psi_4+\frac{1}{5}\Psi_3\Psi_5-\frac{1}{35}\Psi_6\br)\br\},\]
with
\[\Psi_n:=\frac{\psi_s^{(n)}}{\psi_s''}\quad (n\geq 3),\qquad F_n:=\frac{f_s^{(n)}}{f_s}\quad(n\geq 1), \quad f=w^{-k-3/2}\]
being derivatives evaluated at $w=w_s$.
Using the values $\Psi_3=-2i$, $\Psi_4=-6$, $\Psi_5=24i$ and $\Psi_6=120$, we find the values
\[B_1(0)=-\frac{11}{24},\quad B_1(1)=-\frac{47}{24},\qquad B_2(0)=\frac{265}{1152}.\]

Combination of the expansions (\ref{e22}) and (\ref{e23}), together with the conjugate form of expansion for $\Upsilon_k^-(n)$ (with saddle point at $w=+i$), then leads to
\[{\cal J}_1\sim \frac{\Lambda(n)}{\sqrt{2\pi n}}\bl\{e^{i\Phi}\sum_{k\geq 0}\frac{A_k}{n^k}\sum_{r\geq0}\frac{B_r(k)}{n^r}+e^{-i\Phi}\sum_{k\geq 0}\frac{{\bar A}_k}{n^k}\sum_{r\geq0}\frac{B_r(k)}{n^r}\br\}\]
\bee\label{e24}
\hspace{0.6cm}=\sqrt{\frac{2}{\pi n}}\,\Lambda(n)\bl\{\cos \Phi \sum_{k\geq0}\frac{C_k}{n^k}+\sin \Phi \sum_{k\geq1}\frac{D_k}{n^k}\br\}\qquad(n\to\infty),
\ee
where 
\[\Phi:=\fs\pi n+a-\fs\pi\nu,\qquad \Lambda(n):=\frac{n!e^n n^{-n-1/2}}{\sqrt{2\pi}}=1+\frac{1}{12n}+\frac{1}{288n^2}+O(n^{-3})\]
by Stirling's formula for $n!$.
The first few coefficients $C_k$ and $D_k$ are given by
\[C_0=1,\quad C_1=B_1(0)-c_1(\nu),\quad C_2=B_2(0)-c_1(\nu) B_1(1)+c_2(\nu)-\f{3}{8}a^2,\]
\[D_1=\fs a.\quad D_2=\fs a(B_1(1)-3c_1(\nu)).\]

Now we deal with the branch-point contribution. First, the contribution from the small circular path of radius $\rho$ about the branch point $u=-a$ is proportional to $\rho^{\nu+1}$, which vanishes as $\rho\to0$ since $\nu\ge0$. The contribution ${\cal J}_2$ from the upper and lower sides of the branch cut is
\bee\label{e24a}
{\cal J}_2=(-)^{n-1}\frac{n!}{2\pi i}\,\int_0^n \frac{J_\nu(xe^{\pi i})-J_\nu(xe^{-\pi i})}{(a+x)^{n+1}}\,dx=(-)^{n-1}n!\frac{\sin \pi\nu}{\pi} \int_0^n \frac{J_\nu(x)}{(a+x)^{n+1}}\,dx,
\ee
since $J_\nu(xe^{\pm\pi i})=e^{\pm\pi i\nu} J_\nu(x)$. Now \cite[6.563]{GR}
\[\int_0^n\frac{J_\nu(x)}{(a+x)^{n+1}}\,dx=\bl\{\int_0^\infty - \int_n^\infty\br\}\frac{J_\nu(x)}{(a+x)^{n+1}}\,dx\]
\[=(-)^{n-1}\frac{\pi a^{-n}}{\sin \pi\nu}(S_1-S_2)+O(n^{-n-1}),\]
where
\[S_1=(-)^{n-1}\frac{\sin \pi\nu}{\pi}\sum_{m\geq0}\frac{(-)^m (\fs a)^{\nu+2m}\g(1+\nu+2m)}{m! \,\g(1+\nu+m)}\,\frac{\g(n-\nu-2m)}{\g(n+1)},\]
\[S_2=\frac{1}{n! \sqrt{\pi}}\sum_{m\geq0}\frac{a^{n+m}}{m!}\,\frac{\g(\frac{n}{2}+\frac{m+1}{2})\g(\frac{n}{2}+\frac{m+2}{2})}{\g(\frac{n}{2}+\frac{m+2+\nu}{2})\g(\frac{n}{2}+\frac{m+2-\nu}{2})}\,\cos \fs\pi(\nu\!-\!n\!+\!m)\]
and we have employed the well-known bound $|J_\nu(x)|\leq1$ for $\nu\geq0$, $x\geq0$ to bound the second integral over $[n,\infty)$. Making use of the fact that $\g(n+\alpha)/\g(n+\beta)\sim n^{\alpha-\beta}$ as $n\to\infty$  \cite[(5.11.12)]{DLMF}, we find
\[S_1\sim(-)^{n-1}\frac{\sin \pi\nu}{\pi n}\sum_{m\geq0}\frac{(-)^m(a/(2n))^{\nu+2m}\g(1+\nu+2m)}{m! \g(1+\nu+m)}=O(n^{-1-\nu})\]
and 
\[S_2\sim\frac{1}{n!}\sqrt{\frac{2}{\pi n}}\sum_{m\geq0}\frac{a^{n+M}}{m!}\\cos \fs\pi(\nu\!-\!n\!+\!m)=O(a^n n^{-1/2}/n!)\]
as $n\to\infty$. Hence, it follows that\footnote{Use of $|J_\nu(x)|\leq1$ for $\nu\geq0$, $x\geq0$ shows immediately that ${\cal J}_2$ in (\ref{e24a}) is {\it bounded} by $n!a^{-n}/n$.} for $\nu\geq0$
\bee\label{e25}
{\cal J}_2=O\bl(\frac{n! a^{-n}}{n^{1+\nu}}\br)\qquad (n\to\infty).
\ee

We remark that, when $\nu$ is an integer, the branch-point contribution is not present. Thus, collecting together the results in (\ref{e24}) and (\ref{e25}) we have the following theorem:
\newtheorem{theorem}{Theorem}
\begin{theorem}$\!\!\!.$\ \ Let $\nu\geq0$, $a>0$ and $n$ be a positive integer. Then
\bee\label{e26}
J_\nu^{(n)}(a)=\left\{\begin{array}{ll} O\bl(\dfrac{n! a^{-n}}{n^{1+\nu}}\br) & \nu\neq 0, 1, 2, \ldots\\
\\
O(n^{-1/2} \cos \Phi) & \nu=0, 1, 2, \ldots\end{array}\right.
\ee
as $n\to\infty$, where $\Phi=\fs\pi n+a-\fs\pi\nu$.
\end{theorem}
The second result in (\ref{e26}) when $\nu=0$ and $a=1$ was given in \cite[p.~126]{PHad}.

In Table 1 we present some numerical values of the absolute error in the computation of $J_\nu^{(n)}(a)$ compared with the asymptotic expansion in (\ref{e24}). The numerical value of the $n$th derivative follows from the expression \cite[(10.6.7)]{DLMF}
\bee\label{e27}
J_\nu^{(n)}(a)=2^{-n}\sum_{r=0}^n(-)^r\bl(\!\!\begin{array}{c}n\\r\end{array}\!\!\br) J_{\nu-n+2r}(a),
\ee
which can be further simplified when $\nu=m$ is an integer since $J_{-m}(a)=(-)^m J_m(a)$.
\begin{table}[th]
\caption{\footnotesize{Values of the absolute error in the computation of $J_\nu^{(n)}(a)$ from (\ref{e24}) for different $n$ when $\nu=0$, $a=1$ and truncation index $k\leq2$.}}\label{t21}
\begin{center}
\begin{tabular}{r||l|l|l|l}
\hline
\mcol{1}{c||}{$k$} & \mcol{1}{c|}{$n=40$} & \mcol{1}{c|}{$n=60$} & \mcol{1}{c|}{$n=100$} & \mcol{1}{c}{$n=200$}\\ 
\hline
&&&&\\[-0.3cm]
0 & $6.968\times 10^{-4}$ & $3.903\times 10^{-4}$ & $1.856\times 10^{-4}$ & $6.674\times 10^{-5}$\\
1 & $6.373\times 10^{-5}$ & $2.341\times 10^{-5}$ & $6.594\times 10^{-6}$ & $1.174\times 10^{-5}$\\
2 & $2.546\times 10^{-6}$ & $6.193\times 10^{-7}$ & $1.041\times 10^{-7}$ & $9.224\times 10^{-9}$\\
[.2cm]\hline\end{tabular}
\end{center}
\end{table}

\vspace{0.6cm}

\begin{center}
{\bf 3. \ The derivatives $K_\nu^{(n)}(a)$}
\end{center}
\setcounter{section}{3}
\setcounter{equation}{0}
\renewcommand{\theequation}{\arabic{section}.\arabic{equation}}
Proceeding in the same manner as in Section 2, we have the $n$th derivative of the modified Bessel function given by
\bee\label{e30}
K_\nu^{(n)}(a)=\frac{n!}{2\pi i}\oint \frac{K_\nu(a+u)}{u^{n+1}}\,du \qquad (0\leq\nu<1),
\ee
where the closed loop surrounds the origin in the positive sense excluding the branch point at $u=-a$. Note that the branch point is present when $\nu=0$ on account of the logarithmic nature of $K_0(x)$ near $x=0$.

We now expand the loop in (\ref{e30}) into a large circular contour of radius $R$ together with the indentation around the branch point $u=-a$. Since $K_\nu(x)\sim e^{-x}/\sqrt{2\pi x}$ for $|x|\to\infty$ in $|\arg\,x|<\pi$, it is readily seen that the above integral possesses a saddle point at $u=-n$ on the branch cut. Choose $R=n$; the saddle-point contribution is then of order
\[\frac{n! e^n n^{-n-1/2}}{\sqrt{2\pi n}}=O(n^{-1/2})\qquad (n\to\infty).\]
The contribution from the small circular path of radius $\rho$ round the branch point is controlled by $\rho^{1-\nu}\to0$ as $\rho\to0$ for $0<\nu<1$; in the case $\nu=0$, the contribution is controlled by $\rho \log\,\rho\to 0$ as $\rho\to0$.

The contribution from the upper and lower sides of the branch cut is
\[{\cal K}=(-)^{n-1} \frac{n!}{2\pi i}\int_0^n\frac{K_\nu(xe^{-\pi i})-K_\nu(xe^{\pi i})}{(a+x)^{n+1}}\,dx.\]
From the result 
$K_\nu(xe^{\pm \pi i})=e^{\mp\pi i\nu} K_\nu(x)\mp\pi i I_\nu(x)$ \cite[(10.34.2)]{DLMF},
we have
\bee\label{e32}
{\cal K}=n!\int_0^n\bl\{\frac{\sin \pi\nu}{\pi}\,K_\nu(x)+I_\nu(x)\br\}\,\frac{dx}{(a+x)^{n+1}}.
\ee
In the appendix it is established that
\[\int_0^n\!\!\frac{K_\nu(x)}{(a+x)^{n+1}}\,dx=O\bl(\frac{a^{-n}}{n^{1-\nu}}\br),\qquad \int_0^n\!\!\frac{I_\nu(x)}{(a+x)^{n+1}}\,dx=O\bl(\frac{a^{-n}}{n^{1+\nu}}\br)\]
as $n\to\infty$. Hence
\[{\cal K}=\frac{\sin \pi\nu}{\pi}\,O\bl(\frac{n!a^{-n}}{n^{1-\nu}}\br)+O\bl(\frac{n!a^{-n}}{n^{1+\nu}}\br)
=O\bl(\frac{n!a^{-n}}{n^{1-\nu}}\br)\qquad (0\leq\nu<1)\]
to yield the following theorem:
\begin{theorem}$\!\!\!.$\ \ Let $0\leq\nu<1$, $a>0$ and $n$ be a positive integer. Then
\bee\label{e31}
K_\nu^{(n)}(a)=O\bl(\frac{n!a^{-n}}{n^{1-\nu}}\br)
\ee
as $n\to\infty$.
\end{theorem}

We remark that the derivatives of the $K$-Bessel function can be evaluated by means of the formula
\cite[(10.29.5)]{DLMF}
\bee\label{e34}
K_\nu^{(n)}(a)=(-2)^{-n} \sum_{r=0}^n \bl(\!\!\begin{array}{c}n\\r\end{array}\!\!\br) K_{\nu-n+2r}(a).
\ee
\vspace{0.6cm}

\begin{center}
{\bf 4. \ The asymptotic expansion of the integrals $H_J(a,x)$ and $H_K(a,x)$}
\end{center}
\setcounter{section}{4}
\setcounter{equation}{0}
\renewcommand{\theequation}{\arabic{section}.\arabic{equation}}
We first consider the integral $H_J(a,x)$ which can be written as
\[H_J(a,x)=\int_0^a e^{-xt}J_\nu(t)\,dt=\bl(\int_0^\infty-\int_a^\infty \br) e^{-xt} J_\nu(t)\,dt,\]
where $a>0$ and we take $\nu\geq0$. Thus  \cite[p.~386(8)]{W}
\bee\label{e41}
H_J(a,x)=\frac{(\sqrt{1+x^2}-x)^\nu}{\sqrt{1+x^2}}-e^{-ax}\int_0^\infty e^{-au} J_\nu(a+u)\,du,
\ee
where we have put $t=a+u$.

The Taylor series expansion of $J_\nu(a+u)$ is
\[J_\nu(a+u)=\sum_{n\geq 0}J_\nu^{(n)}(a)\,\frac{u^n}{n!},\]
where, from Theorem 1 describing the behaviour of $J_\nu^{(n)}(a)$ for large $n$, it is seen that the series converges for $|u|<a$ when $\nu>0$ is non-integer, but converges for $|u|<\infty$ when $\nu=0, 1, 2, \ldots\ $. When $\nu>0$ is non-integer the integral on the right-hand side of (\ref{e41}) becomes
\[e^{-ax}\bl\{\sum_{n\geq0}\frac{J_\nu^{(n)}(a)}{n!} \int_0^au^n e^{-xu}du+T_J(a,x)\br\}\]
\bee\label{e42}
=e^{-ax}\bl\{\sum_{n\geq0}\frac{J_\nu^{(n)}(a)}{x^{n+1}}\,P(n+1,ax)+T_J(a,x)\br\},
\ee
where $P(\alpha,x)=\gamma(\alpha,x)/\g(\alpha)$ is the normalised incomplete gamma function and
\[T_J(a,x):=\int_a^\infty e^{-xu} J_\nu(a+u)\,du.\]
Use of the bound $|J_\nu(x)|<1$  $ (\nu>0,\ x>0)$ shows that $|T_J(a,x)|<e^{-ax}/x$.

The series in (\ref{e42}) is an example of a Hadamard expansion; see \cite[Ch.~2]{PHad} for a full discussion of the use of such expansions in hyperasymptotic evaluation. The presence of the incomplete gamma function in this series acts as a `smoothing' factor on the coefficients $J_\nu^{(n)}(a)/x^{n+1}$, since the behaviour of $P(\alpha,x)$ is given by
\[P(\alpha,x)\sim\left\{\begin{array}{ll} 1 & (x\to\infty)\\
\dfrac{x^\alpha e^{-x}}{\g(1+\alpha)} & (\alpha\to\infty).\end{array}\right.\]
Thus, $P(\alpha,x)$ changes from approximately unity when $\alpha\,\ltwid\, x$ to a rapid decay to zero when $\alpha \,\gtwid\,x$. Consequently, the early terms ($n<ax$) in the above series behave like those of the associated Poincar\'e asymptotic series
\[e^{-ax}\sum_{n\geq 0} \frac{J_\nu^{(n)}(a)}{x^{n+1}}\qquad (x\to\infty).\]
Hence, for $\nu>0$ and non-integer we have the asymptotic expansion
\bee\label{e43}
H_J(a,x)-\frac{(\sqrt{1+x^2}-x)^\nu}{\sqrt{1+x^2}}\sim -e^{-ax}\sum_{n\geq 0} \frac{J_\nu^{(n)}(a)}{x^{n+1}}
\ee
as $x\to\infty$ and fixed finite $a$. 

For integer values of $\nu$ we have the {\it absolutely convergent} expansion
\bee\label{e44}
H_J(a,x)=\frac{(\sqrt{1+x^2}-x)^\nu}{\sqrt{1+x^2}}-e^{-ax}\sum_{n\geq0} \frac{J_\nu^{(n)}(a)}{x^{n+1}}\qquad (ax>1;\ \nu=0, 1, 2, \ldots).
\ee
The case of (\ref{e44}) when $\nu=0$, $a=1$ has been given\footnote{There is a misprint in \cite[(2.3.5)]{PHad}: the $x^{-n}$ should be $x^{-n-1}$.} in \cite[p.~125]{PHad}.

A similar treatment of the integral $H_K(a,x)$ yields, when $0\leq\nu<1$,
\bee\label{e45}
H_K(a,x)=\int_0^a e^{-xt} K_\nu(t)\,dt=H_K(\infty,x)-e^{-ax}\int_0^\infty e^{-xu} K_\nu(a+u)\,du,
\ee
where \cite[p.~388(9)]{W}
\bee\label{e46}
H_K(\infty,x)=\frac{\pi}{2\sin \pi\nu}\,\frac{1}{\sqrt{x^2-1}}\bl\{(x+\sqrt{x^2-1})^\nu-(x+\sqrt{x^2-1})^{-\nu}\br\}.
\ee
The Taylor series series expansion of $K_\nu(a+u)$ given by
\[K_\nu(a+u)=\sum_{n\geq0} \frac{K_\nu^{(n)}(a)}{n!}\,u^n\]
converges for $|u|<a$ by Theorem 2. Thus, the integral appearing on the right-hand side of (\ref{e45}) becomes
\[e^{-ax}\bl\{\sum_{n\geq0}\frac{K_\nu^{(n)}(a)}{x^{n+1}}\,P(n+1,ax)+T_K(a,x)\br\},\]
where, by the inequality in (\ref{a0}) valid for $\nu>0$,
\[T_K(a,x):=\int_a^\infty e^{-xu} K_\nu(a+u)\,du < 2^{\nu-1} \g(\nu)\int_a^\infty e^{-xu} (a+u)^{-\nu}du \]
\[< 2^{\nu-1} \g(\nu)\int_a^\infty e^{-xu} u^{-\nu}du =\frac{2^{\nu-1} \g(\nu)}{x^{1-\nu}}\,\g(1-\nu,ax).\]
From \cite[(8.11.2)]{DLMF} the upper incomplete gamma function has the behaviour $\g(1-\nu,ax)\sim (ax)^{-\nu} e^{-ax}$ as $ax\to\infty$, so that $T_K(a,x)$ is exponentially small\footnote{When $\nu=0$, the bound $K_0(x)<(2/\pi x)^{1/2}e^{-x}$ ($x>0$)
 \cite{REG} shows that $T_K(a,x)$ is also exponentially small in $x$ in this case.} as $x\to\infty$ with $a$ bounded away from zero.

Hence we obtain the asymptotic expansion
\bee\label{e47}
H_K(a,x)-H_K(\infty,x)
\sim -e^{-ax}\sum_{n\geq0}\frac{K_\nu^{(n)}(a)}{x^{n+1}}
\ee
as $x\to\infty$ with fixed finite $a$, where $H_K(\infty,x)$ is defined in (\ref{e46}).
\vspace{0.6cm}

\begin{center}
{\bf 5. \ Concluding remarks}
\end{center}
\setcounter{section}{5}
\setcounter{equation}{0}
\renewcommand{\theequation}{\arabic{section}.\arabic{equation}}
We have investigated the high-order derivatives of the Bessel functions $J_\nu(a)$ and $K_\nu(a)$ for $a>0$ and applied these results to determine the asymptotic character of the expansion of two incomplete Laplace transforms of these functions as the transform variable $x\to+\infty$. In the case of the integer order $J$ Bessel function the expansion was shown to be convergent. The connection of these expansions to the recently developed theory of Hadamard
expansions has been indicated. However, although hyperasymptotic precision is possible with this latter procedure, this aspect is not pursued here.  

A similar procedure can be employed to determine the expansion of the integrals involving the Bessel functions $Y_\nu(t)$ and $I_\nu(t)$ given by
\[H_Y(a,x)=\int_0^ae^{-xt} Y_\nu(t)\,dt \quad(-1<\nu<1),\qquad H_I(a,x)=\int_0^a e^{-xt}I_\nu(t)\,dt \quad(x>1,\ \nu\geq0).\]
Then we have the expansions as $x\to+\infty$
\bee
H_Y(a,x)-H_Y(\infty,x)\sim -e^{-ax}\sum_{n\geq0}\frac{Y_\nu^{(n)}(a)}{x^{n+1}}\qquad(-1<\nu<1),
\ee
where
\[H_Y(\infty,x)=-\frac{\csc \pi\nu}{\sqrt{1+x^2}}\,\bl\{(x+\sqrt{1+x^2})^{\nu}-\cos \pi\nu\,(x+\sqrt{1+x^2})^{-\nu}\br\},\]
and
\bee
H_I(a,x)-H_I(\infty,x)\sim -e^{-ax}\sum_{n\geq0}\frac{I_\nu^{(n)}(a)}{x^{n+1}}\qquad(\nu\neq 0, 1, 2, \ldots),
\ee
where
\[H_I(\infty,x)=\frac{(x+\sqrt{x^2-1})^{-\nu}}{\sqrt{x^2-1}}\qquad (x>1).\]
The derivatives $Y_\nu^{(n)}(a)$ and $I_\nu^{(n)}(a)$ can be obtained from (\ref{e27}) with $J_\nu$ replaced by $Y_\nu$, and from (\ref{e34}) with $K_\nu$ replaced by $I_\nu$ and the factor $(-2)^{-n}$ replaced by $2^{-n}$.

When $\nu=m$, where $m$ is a non-negative integer, we have the convergent expansion
\bee
H_I(a,x)=H_I(\infty,x)-e^{-ax}\sum_{n\geq0}\frac{I_m^{(n)}(a)}{x^{n+1}}\qquad (x>1),
\ee
since routine calculations similar to those described in Section 2 show that
\[I_m^{(n)}(a) \sim \frac{1}{\sqrt{2\pi n}}\,(e^a+(-)^{m+n}e^{-a})\qquad (n\to\infty).\]

\vspace{0.6cm}

\begin{center}
{\bf Appendix: Estimation of two integrals appearing in ${\cal K}$}
\end{center}
\setcounter{section}{1}
\setcounter{equation}{0}
\renewcommand{\theequation}{\Alph{section}.\arabic{equation}}
In this appendix we estimate the growth of the integrals appearing in ${\cal K}$ efined in (\ref{e31}) for large $n$. 
Consider first the integral, where $a>0$ is a fixed parameter,
\[{\cal K}_1=\int_0^n\frac{K_\nu(x)}{(a+x)^{n+1}}\,dx=a^{-n}\int_0^{n/a}\frac{K_\nu(ax)}{(1+x)^{n+1}}\,dx.\]
From \cite{REG}, we have the bounds
\bee\label{a0}
2^{\nu-1}\g(\nu)e^{-x}<x^\nu K_\nu(x)<2^{\nu-1}\g(\nu)\qquad (\nu>0,\ x>0),
\ee
whence, for $0<\nu<1$,
\bee\label{a1}
\frac{2^{\nu-1}\g(\nu)}{a^{n+\nu}}\int_0^{n/a}\!\!\!\frac{x^{-\nu} e^{-ax}}{(1+x)^{n+1}}\,dx=L_1<{\cal K}_1<U_1=\frac{2^{\nu-1}\g(\nu)}{a^{n+\nu}}\int_0^{n/a} \!\!\!\frac{x^{-\nu}}{(1+x)^{n+1}}\,dx.
\ee

Now, since $(1+x)^{-\alpha}>e^{-\alpha x}$ ($\alpha>0$), we have
\[\int_0^{n/a} \frac{x^{-\nu}e^{-ax}}{(1+x)^{n+1}}\,dx>\int_0^{n/a}x^{-\nu} e^{-(n+1+a)x}dx=\frac{1}{(n+1+a)^{1-\nu}}\int_0^X \tau^{-\nu}e^{-\tau}d\tau,\]
where $X=n(n+1+a)/a$. Evaluation of this last integral as the lower incomplete gamma function $\gamma(1-\nu,X)$ and use of its asymptotic behaviour for large $X$ \cite[(8.11.2)]{DLMF} shows that as $n\to\infty$
\[\int_0^{n/a} \frac{x^{-\nu}e^{-ax}}{(1+x)^{n+1}}\,dx>\frac{\g(1-\nu)}{(n+1+a)^{1-\nu}}\{1-O(n^{-2\nu}e^{-n^2})\}.\]
Hence the lower bound satisfies
\bee\label{a2}
L_1>2^{\nu-1} \frac{\pi}{\sin \pi\nu}\,\frac{a^{-n-\nu}}{(n+1+a)^{1-\nu}}=O\bl(\frac{a^{-n}}{n^{1-\nu}}\br)\qquad (n\to\infty).
\ee

For the upper bound $U_1$ we have
\[\int_0^{n/a}\frac{x^{-\nu}}{(1+x)^{n+1}}\,dx=\g(1-\nu)\,\frac{\g(n+\nu)}{\g(n+1)}-\int_{n/a}^\infty\frac{x^{-\nu}}{(1+x)^{n+1}}\,dx\]
\[=\frac{\g(1-\nu)}{n^{1-\nu}}(1+O(n^{-1}))-O(n^{-n-1-\nu})\]
as $n\to\infty$. This yields the extimate
\bee\label{a3}
U_1\sim 2^{\nu-1}\frac{\pi}{\sin \pi\nu}\,\frac{a^{-n-\nu}}{n^{1-\nu}}.
\ee
Consequently, from (\ref{a1})--(\ref{a3}) it follows that
\bee\label{a4}
\int_0^n\frac{K_\nu(x)}{(a+x)^{n+1}}\,dx=O\bl(\frac{a^{-n}}{n^{1-\nu}}\br)\qquad (n\to\infty).
\ee

For the integral
\[{\cal K}_2=\int_0^n\!\!\frac{I_\nu(x)}{(1+x)^{n+1}}\,dx=a^{-n}\int_0^{n/a}\!\!\frac{I_\nu(ax)}{(1+x)^{n+1}}\,dx,\]
we employ the bounds \cite{L}
\[\frac{2^{-\nu}}{\g(1+\nu)}\leq x^{-\nu}I_\nu(x)<\frac{2^{-\nu} e^x}{\g(1+\nu)}\qquad (\nu\geq0,\ x\geq0).\]
Then we obtain
\bee\label{a5}
\frac{2^{-\nu}a^{-n+\nu}}{\g(1+\nu)}\int_0^{n/a}\!\!\!\frac{x^\nu}{(1+x)^{n+1}}\,dx=L_2<{\cal K}_2<U_2=
\frac{2^{-\nu}a^{-n+\nu}}{\g(1+\nu)}\int_0^{n/a}\!\!\!\frac{x^\nu e^{ax}}{(1+x)^{n+1}}\,dx.
\ee

The integrand of the integral on the right-hand side of (\ref{a5}) has, when $\nu>0$, a maximum at $x\simeq \nu/n$ and an absolute minimum at $x\simeq n/a$ (corresponding to the saddle point in (\ref{e40})). Beyond this minimum point the integrand thereafter steadily increases. To estimate this integral for large $n$, we divide the integration path into
$[0,n^{-\mu}]$ and $[n^{-\mu},n/a]$, where $\fs<\mu<1$. The value of the integrand at $x=n^{-\mu}$ is O($n^{-\mu\nu}\exp\,[an^{-\mu}-n^{1-\mu}])$ for large $n$. The contribution from the interval $[n^{-\mu},n/a]$ is therefore O($n^{1-\mu\nu} \exp\,[an^{-\mu}-n^{1-\mu}])$, which is seen to be exponentially small as $n\to\infty$.

Now
\[\frac{x^\nu e^{ax}}{(1+x)^{n+1}}=x^\nu e^{-(n+1-a)x}\,e^{(n+1)(\fs x^2-p(x))}>x^\nu e^{-(n+1-a)x} e^{(n+1)x^2/2},\]
since $p(x):=\log (1+x)-x+\fs x^2>0$ for $x>0$. Then
\[\frac{x^\nu e^{ax}}{(1+x)^{n+1}}>x^\nu e^{-(n+1-a)x}\{1+O(n^{1-2\mu})\}, \qquad x\in[0,n^{-\mu}], \]
and hence the contribution to ${\cal K}_2$ from the interval $[0,n^{-\mu}]$ is
\[\int_0^{n^{-\mu}}\!\!\!\frac{x^\nu e^{ax}}{(1+x)^{n+1}}\,dx>\int_0^{n^{-\mu}} \!\!\!x^\nu e^{-(n+1-a)x}\{1+O(n^{1-2\mu})\}\,dx\]
\[\hspace{1.8cm}=\frac{\gamma(1+\nu,X)}{(n\!+\!1\!-\!a)^{1+\nu}}\{1+O(n^{1-2\mu})\},\]
where $X:=n^{-\mu}(n\!+\!\!1\!-a)\sim n^{1-\mu}$.
Thus, $X\to\infty$ as $n\to\infty$, and we find
\bee\label{a6}
\int_0^{n^{-\mu}}\!\!\!\frac{x^\nu e^{ax}}{(1+x)^{n+1}}\,dx
=\frac{\g(1+\nu)}{(n\!+\!1\!-\!a)^{1+\nu}}+O\bl(\frac{1}{n^{2\mu+\nu}}\br)=O\bl(\frac{1}{n^{1+\nu}}\br)
\ee
with the correction terms from $\gamma(1+\nu,X)$ being exponentially small as $n\to\infty$. From (\ref{a5}) and (\ref{a6}) (applied when $a>0$ and $a=0$), it then follows that
\[\frac{2^{-\nu}a^{-n+\nu}}{(n\!+\!1)^{1+\nu}}<{\cal K}_2<\frac{2^{-\nu}a^{-n+\nu}}{(n\!+\!1\!-\!a)^{1+\nu}},\]
whence
\bee\label{a7}
\int_0^n\!\!\frac{I_\nu(x)}{(1+x)^{n+1}}\,dx=O\bl(\frac{a^{-n}}{n^{1+\nu}}\br)\qquad (n\to\infty).
\ee

\end{document}